\documentclass[12pt]{amsart}
\usepackage{amsmath}
\usepackage{amsfonts}

\theoremstyle{plain}

\newtheorem{corollary}{Corollary}

\newtheorem{proposition}{Proposition}
\newtheorem{remark}{Remark}

\numberwithin{equation}{section}
\input{tcilatex}

\begin{document}
\title[On Hyper Singular Integral Operators]{On Hyper Singular Integral
Operators over Weighted Sobolev Spaces}
\author{Dejenie A. Lakew}
\address{Virginia Union University \\
Department of Mathematics \& Computer Science\\
Richmond, Virginia 23220}
\email{dalakew@vuu.edu}
\urladdr{http://www.vuu.edu}
\date{August 8, 2009}
\dedicatory{}
\thanks{This paper is in final form and no version of it will be submitted
for publication elsewhere.}

\begin{abstract}
In this paper we study singular integral operators which are hyper or weak
over Lipscitz or H\"{o}lder spaces and over weighted Sobolev spaces defined
on unbounded domains in the standard $n$-$D$ Euclidean space $%
\mathbb{R}
^{n}$ for $n\geq 1$. The $\pi -$operator in this case is one of the hyper
singular integral operators which has been studied extensibly than other
hyper singular integral operators. It will be shown the control of
singularity such integral operators that are defined interms of Cauchy
generating kernels by working on weighted Sobolev spaces $W^{p,k}\left(
\Omega ,\parallel x\parallel ^{\zeta +\epsilon }dx\right) $ for some $%
\epsilon >0$ and $\zeta $ some positive integer.
\end{abstract}

\maketitle

\section{Hyper\textbf{\ Singular Integral Operators}}

In this short note we discuss few points about super singular integral
operators, weak(or sub) singular and just singular integral operators by
showing few examples and present some results.

\ \ 

We therefore introduce general singular integral operators in terms of
integrals with Cauchy generating kernels and some other general singular
integral operators with out kernels.

\ 

The calculus versions of singular integral operators are improper integrals,
integrals with unbounded integrands or integrals with unbounded intervals of
integrations.

\ 

To start our work, let $\Omega $ be some bounded domain in the Euclidean
space $%
\mathbb{R}
^{n}$ and $\psi $ be some integrable function over $\Omega $ and $x_{0}\in
\Omega ^{\text{int}}$, interior of the domain with the property that 
\begin{equation*}
\underset{x\rightarrow x_{0}}{\lim }\mid \psi \left( x\right) \mid =\infty
\end{equation*}%
which in this case $x_{0}$ is a singular point of the function.

\ \ \ 

The integral given by 
\begin{equation*}
\dint\limits_{\Omega }\psi \left( x\right) dx
\end{equation*}%
is called a singular integral of the function $\psi $ over the domain $%
\Omega $ with a singularity point $x_{0}$. We evaluate such singular
integrals by evaluating the Cauchy principal value of the singular integral
which is given as follows.

\ \ \ 

Let $\epsilon >0$ and consider the ball $B\left( x_{0},\epsilon \right) $
and define $\Omega _{\epsilon }:=\Omega \backslash B\left( x_{0},\epsilon
\right) $. Then we consider the integral over the deleted sub domain $\Omega
_{\epsilon }$ by 
\begin{equation*}
\dint\limits_{\Omega _{\epsilon }}\psi \left( x\right) dx
\end{equation*}
which avoids the singularity $x_{0}.$\ 

\ \ 

If the limit : 
\begin{equation*}
\underset{\epsilon \rightarrow 0}{\lim }\dint\limits_{\Omega _{\epsilon
}}\psi \left( x\right) dx
\end{equation*}%
called the Cauchy principal value(c.p.v.) exits, then we define the value of
the singular integral as: 
\begin{equation*}
\dint\limits_{\Omega }\psi \left( x\right) dx:=\underset{\epsilon
\rightarrow 0}{\lim }\dint\limits_{\Omega _{\epsilon }}\psi \left( x\right)
dx
\end{equation*}

Examples \ of elementary singular integral operators are given below:

\ \ 

In the unidimensional Euclidean space $%
\mathbb{R}
^{1}:$ let $\Omega =\left( -1,1\right) $ and define the function by 
\begin{equation*}
\psi _{\alpha }\left( x\right) =\mid x\mid ^{-\alpha },\text{ \ for }\
0<\alpha <1
\end{equation*}

Then the function $\psi _{\alpha }$ has a singularity at $0$, since 
\begin{equation*}
\underset{x\rightarrow 0}{\lim }\mid \psi _{\alpha }\left( x\right) \mid
=\infty
\end{equation*}

Therefore, the integral given by 
\begin{equation*}
\dint\limits_{\Omega }\psi _{\alpha }\left( x\right) dx
\end{equation*}
is a weakly singular integral

Let $\epsilon >0$ and consider 
\begin{equation*}
\Omega _{\epsilon }=\Omega \backslash B\left( 0,\epsilon \right) =\left(
-1,1\right) \backslash \left( -\epsilon ,\epsilon \right) .
\end{equation*}

Then the integral $\dint\limits_{\Omega _{\epsilon }}\psi _{\alpha }\left(
x\right) dx$ is no more a singular integral at $x_{0}$ and therefore has a
finite integral as long as the function $\psi $ is integrable on the domain $%
\Omega $.

\ \ 

Therefore, 
\begin{equation*}
\dint\limits_{\Omega _{\epsilon }}\psi _{\alpha }\left( x\right)
dx=\dint\limits_{\Omega _{\epsilon }}\mid x\mid ^{-\alpha }dx
\end{equation*}

is a function of $\alpha $ and $\epsilon $ and if we denote this function by 
$I\left( \alpha ,\epsilon \right) $, then we have 
\begin{equation*}
I\left( \alpha ,\epsilon \right) =\frac{2}{1-\alpha }\left( 1-\epsilon
^{1-\alpha }\right)
\end{equation*}%
which is a finite value in terms of $\epsilon $ and $\alpha $. Then taking
the c.p.v. of the above integral :

\begin{eqnarray*}
\underset{\epsilon \rightarrow 0}{\lim }\dint\limits_{\Omega _{\epsilon
}}\psi _{\alpha }\left( x\right) dx &=&\underset{\epsilon \rightarrow 0}{%
\lim }\dint\limits_{\Omega _{\epsilon }}\mid x\mid ^{-\alpha }dx \\
&=&\underset{\epsilon \rightarrow 0}{\lim }I\left( \alpha ,\epsilon \right)
\\
&=&\frac{2}{1-\alpha }
\end{eqnarray*}%
as $1-\alpha >0$.

When $\alpha =1,$ the function is $\psi _{-1}\left( x\right) =\mid x\mid
^{-1}$ and this function generates an integral 
\begin{equation*}
\dint\limits_{\Omega }\psi _{-1}\left( x\right) dx
\end{equation*}
called a singular integral.

For $\alpha =1+\varepsilon $, $\varepsilon >0$, the integral 
\begin{equation*}
\dint\limits_{\Omega }\psi _{\alpha }\left( x\right) dx
\end{equation*}
is called a hyper singular integral. Besides 
\begin{eqnarray*}
\underset{\epsilon \rightarrow 0}{\lim }\dint\limits_{\Omega _{\epsilon
}}\psi _{\alpha }\left( x\right) dx &=&\underset{\epsilon \rightarrow 0}{%
\lim }\dint\limits_{\Omega _{\epsilon }}\mid x\mid ^{-\alpha }dx \\
&=&\underset{\epsilon \rightarrow 0}{\lim }I\left( \alpha ,\epsilon \right)
\\
&=&\underset{\epsilon \rightarrow 0}{\lim }\frac{2}{1-\alpha }\left(
1-\epsilon ^{1-\alpha }\right) =\infty
\end{eqnarray*}

\bigskip Therefore the improper integral is divergent.

We therefore construct the classical singular integral operators which are
obtained from generating Kernels.

\ \ 

Let us begin with one of the most common generating kernels given by the
function:%
\begin{equation*}
K\left( x\right) =\frac{\overline{x}}{\omega _{n}\mid x\mid ^{n}}
\end{equation*}%
which is called the Cauchy kernel whose singularity is at zero.

\ \ \ 

This kernel gives singular integral operator on the space of functions such
that the convolution is finite over the domain $\Omega $, which is given by 
\begin{equation*}
\Phi \left( \psi \right) \left( x\right) =\dint\limits_{\Omega }K\left(
x-y\right) \psi \left( y\right) d\Omega _{y}
\end{equation*}

From the classification of singular integrals, we will see that $\Phi $ is
indeed a weak singular integral: let $\lambda \in 
\mathbb{R}
_{>0}$, 
\begin{equation*}
K\left( \lambda x\right) =\frac{\lambda \overline{x}}{\omega _{n}\lambda
^{n}\mid x\mid ^{n}}=\lambda ^{-\left( n-1\right) }K(x)
\end{equation*}%
which gives that $K$ is a homogeneous function of exponent $n-1$ which is
less than $n.$ The singular integral operator $\Phi $ given above in
literature is called the Teodorescu transform.

\ 

It is an important transform in Sobolev spaces with a regularity
augmentation property by one : 
\begin{equation*}
\Phi :W^{p,k}\left( \Omega \right) \rightarrow W^{p,k+1}\left( \Omega
\right) .
\end{equation*}%
We can further study the function spaces where the weak singular integral
works. In the sequel, we use the following set up:

For $\varepsilon >0,$ consider $B\left( x,\varepsilon \right) ,$ the $%
\varepsilon -$ball centered at $x$ and radius $\varepsilon $ and consider
the punctured domain $\Omega _{\varepsilon }=\Omega \backslash B\left(
x,\varepsilon \right) .$

\ \ \ 

\begin{proposition}
If $\Omega $ is unbounded and smooth domain in $%
\mathbb{R}
^{n}$, then $K\left( x\right) $ is $p-$integrable over $\Omega _{\varepsilon
}$ for $\frac{n}{n-1}<p<\infty $.
\end{proposition}

\begin{proof}
\begin{equation*}
\Vert K\left( x\right) \Vert =\Vert \frac{\overline{x}}{\omega _{n}\mid
x\mid ^{n}}\Vert =\frac{r^{1-n}}{\omega _{n}}
\end{equation*}%
for $\Vert x\Vert =r$ and using polar coordinates, we have the following
norm estimates:

\begin{equation*}
\dint\limits_{\Omega _{\varepsilon }}\Vert K\left( x\right) \Vert ^{p}dx\leq
c\left( \theta \right) \dint\limits_{\varepsilon }^{\infty }r^{p\left(
1-n\right) +n-1}dr
\end{equation*}%
\begin{equation*}
=c\left( \theta \right) \underset{\sigma \rightarrow \infty }{\lim }\left( 
\frac{r^{p\left( 1-n\right) +n}}{p(1-n)+n}\mid _{\varepsilon }^{\sigma
}\right)
\end{equation*}

\begin{equation*}
=c(\theta )\underset{\sigma \rightarrow \infty }{\lim }\left( \frac{\sigma
^{p\left( 1-n\right) +n}}{p(1-n)+n}-\frac{\varepsilon ^{p\left( 1-n\right)
+n}}{p(1-n)+n}\right)
\end{equation*}%
and this is finite\ and equals $c(\theta )\left( \frac{\varepsilon ^{p\left(
1-n\right) +n}}{p(n-1)-n}\right) ,$ if 
\begin{equation*}
p(1-n)+n<0
\end{equation*}

That is 
\begin{equation*}
\frac{n}{n-1}<p<\infty
\end{equation*}
which proves the proposition.
\end{proof}

\bigskip\ \ \ \ 

If the domain is a bounded smooth one, then we consider a singularity at a
finite point and the exponent of integrability will be different.

\ \ \ \ 

Now, as we see that $K$ is in the Sobolev space $W^{p,k}\left( \Omega
_{\varepsilon }\right) $ for $p>\frac{n}{n-1}$, we can determine the
function space where we can work with this function as a generating kernel
for singular integral operators .

\ 

\begin{proposition}
The convolution $K\ast \mid _{\Omega _{\varepsilon }}f$ is well defined and
finite over $W^{q,k}\left( \Omega _{\varepsilon }\right) $ for $1<q<n$.
\end{proposition}

\begin{proof}
From H\"{o}lder's inequality, the product $Kf\in W^{1,k}\left( \Omega
_{\varepsilon }\right) $ when $K\in W^{p,k}\left( \Omega _{\varepsilon
}\right) $ and $f\in W^{q,k}\left( \Omega _{\varepsilon }\right) $ such that 
$p^{-1}+q^{-1}=1$.

Therefore as $p\in (\frac{n}{n-1},\infty )$, we have $q\in (1,n)$ which is
the required result.
\end{proof}

\ \ 

\begin{remark}
One can see that\ in $%
\mathbb{R}
^{2}$, $p$ should strictly be greater than $2$ and therefore we can not work
over the function space $W^{2,k}\left( \Omega \right) $ using the kernel as
it is.
\end{remark}

\begin{proposition}
Let $\Omega $ be a smooth, unbounded domain in $%
\mathbb{R}
^{n}$ and $p\in (\frac{n}{n-1},\infty )$ and $q$ be the conjugate index of $%
p $. Then we have :

\begin{equation*}
\parallel Kf\parallel _{W^{1,k}\left( \Omega _{\epsilon }\right) }\leq
\parallel K\parallel _{W^{p,k}\left( \Omega _{\varepsilon }\right)
}\parallel f\parallel _{W^{q,k}\left( \Omega _{\varepsilon }\right)
}\rightarrow \parallel K\parallel _{W^{p,k}\left( \Omega _{\varepsilon
}\right) }\parallel f\parallel _{W^{n,k}\left( \Omega _{\varepsilon }\right)
}
\end{equation*}%
as $q\nearrow n$.
\end{proposition}

\begin{proof}
From H\"{o}lder's inequality, we have 
\begin{equation*}
\parallel Kf\parallel _{W^{1,k}\left( \Omega _{\varepsilon }\right)
}=\dint\limits_{\Omega _{\varepsilon }}\mid Kf\mid \leq \parallel K\parallel
_{W^{p,k}\left( \Omega _{\varepsilon }\right) }.\parallel f\parallel
_{W^{q,k}\left( \Omega _{\varepsilon }\right) }
\end{equation*}%
Then taking the limiting norm on the indices $p$ and $q$ with $%
p^{-1}+q^{-1}=1$ we have :

\begin{equation*}
\underset{q\nearrow n}{\lim }\left( \parallel K\parallel _{W^{p,k}\left(
\Omega _{\varepsilon }\right) }.\parallel f\parallel _{W^{q,k}\left( \Omega
_{\varepsilon }\right) }\right) =\parallel K\parallel _{W^{p,k}\left( \Omega
_{\varepsilon }\right) }.\parallel f\parallel _{W^{n,k}\left( \Omega
_{\varepsilon }\right) }
\end{equation*}%
since $p\searrow \left( \frac{n}{n-1}\right) \Rightarrow q\nearrow n$ and
that finishes the argument.
\end{proof}

\bigskip\ \ \ 

The next singular integral we consider is the one generated from the
fundamental solution of the Laplacian operator%
\begin{equation*}
\Delta =-\dsum\limits_{j=1}^{n}\frac{\partial ^{2}}{\partial x_{j}^{2}}%
=\left( \dsum\limits_{j=1}^{n}e_{j}\frac{\partial }{\partial x_{j}}\right) 
\overline{\left( \dsum_{j=1}^{n}e_{j}\frac{\partial }{\partial x_{j}}\right) 
}
\end{equation*}%
which is given by 
\begin{equation*}
\Psi _{2.\Omega }\left( x\right) =\frac{-1}{\omega _{n}\parallel x\parallel
^{n+2}}
\end{equation*}
and the corresponding singular integral associated is given by 
\begin{equation*}
\Phi _{2,\Omega }\left( \phi \right) =\dint\limits_{\Omega }\Psi _{2.\Omega
}\left( x-y\right) \phi \left( y\right) d\Omega _{y}
\end{equation*}

We investigate in which generalized Lebesgue space is $\Psi _{2.\Omega }$
over unbounded domain $\Omega \subseteq 
\mathbb{R}
^{n}.$

\ \ \ 

\begin{proposition}
Let $\Omega $ be a smooth and unbounded domain in $%
\mathbb{R}
^{n}$ for $n\geq 1$. Then $\Psi _{2.\Omega }\in W^{p,k}\left( \Omega
_{\varepsilon },Cl_{n}\right) $ for $p\in (\frac{n}{n+2},\infty )$.
\end{proposition}

\begin{proof}
Consider the integral 
\begin{equation*}
\dint\limits_{\Omega _{\delta }}\mid \Psi _{2.\Omega }\mid ^{p}dx
\end{equation*}%
using polar coordinates, the integral becomes $:$%
\begin{equation*}
c\left( \theta ,\omega _{n}\right) \dint\limits_{\delta }^{\infty
}r^{-\left( n+2\right) p+n-1}dr
\end{equation*}

\ \ 

and it will be finite towards the boundary of the domain when $p>\frac{n}{n+2%
}$, where $c\left( \theta ,\omega _{n}\right) $ is a constant that depends
on $\theta $ and the surface area $\omega _{n}$ of the unit sphere $S^{n-1}.$
\end{proof}

\bigskip\ \ \ 

\section{\textbf{\ Weighted Sobolev Spaces}}

If we try to find Sobolev spaces in which the kernel $\Psi _{2.\Omega }$\
works, we might end up in working with a dual spaces whose conjugate indices
are negative.

\bigskip\ \ \ 

For instance in the limiting cases : $q\rightarrow \frac{-n}{2}$ as $%
p\searrow \frac{n}{n+2}$, which shows that $q$ has a negative limiting index
which is going to be a conjugate index of a limiting index of $p$ in some
sense.

\ \ 

To remedy this, we introduce a weight on the Lebesgue volume measure $dx$ so
that we avoid dual spaces with negative indices.

\ 

The weight function that we choose stretches the Lebesgue volume measure so
that the singularity from the kernel is better managed and made more
controlled.

\bigskip\ \ \ 

We choose a radial weight function given by $w(x)=\parallel x\parallel
^{2+\varepsilon }$, where $\varepsilon $ is some positive constant and we
investigate the integral : 
\begin{equation*}
\dint\limits_{\Omega _{\delta }}\Psi _{2.\Omega }(x)d\mu (x)
\end{equation*}%
where $d\mu \left( x\right) =w(x)dx$.

\ \ 

\begin{proposition}
Over unbounded domain $\Omega \subseteq 
\mathbb{R}
^{n}$, $\Psi _{2.\Omega }\in W^{p,k}\left( \Omega _{\delta }\right) $ for $1+%
\frac{\varepsilon }{n+2}<p<\infty $.
\end{proposition}

\begin{proof}
We see from the proposition that the interval for the index $p$ is much
improved and the conjugate space will be a dual space with positive index.

\ \ 

Therefore,%
\begin{eqnarray*}
\dint\limits_{\Omega _{\delta }} &\parallel &\Psi _{2.\Omega }\parallel
^{p}d\mu \left( x\right) =c\left( \theta ,\omega _{n}\right)
\dint\limits_{\delta }^{\infty }r^{-\left( n+2\right) p+2+n+\varepsilon -1}dr
\\
&=&c\left( \theta ,\omega _{n}\right) \underset{\rho \rightarrow \infty }{%
\lim }\left( \frac{r^{-\left( n+2\right) p+n+2+\varepsilon }}{-\left(
n+2\right) p+n+2+\varepsilon }\mid _{\delta }^{\rho }\right) <\infty
\end{eqnarray*}%
when $-\left( n+2\right) p+n+2+\varepsilon <0$ which implies that $p>\frac{%
2+n+\varepsilon }{2+n}$ which is the required result.
\end{proof}

\bigskip\ \ \ \ 

Next, we determine the Sobolev space $W^{q,k}\left( \Omega \right) $ in
which the product $\Psi _{2.\Omega }\phi $ is integrable or the convolution $%
\Psi _{2.\Omega }\ast _{\mid \Omega }\phi $ is finite.

\ \ 

\begin{proposition}
Over unbounded domain $\Omega \subseteq 
\mathbb{R}
^{n}$, and for $1+\frac{\varepsilon }{n+2}<p<\infty $, with respect to the
weighted measure $d\mu \left( x\right) =\parallel x\parallel ^{2+\varepsilon
}dx$ we have

\begin{equation*}
\Psi _{2.\Omega }\phi \in W^{1,k}\left( \Omega ,\parallel x\parallel
^{2+\varepsilon }dx\right)
\end{equation*}
when 
\begin{equation*}
\phi \in W^{q,k}\left( \Omega ,\parallel x\parallel ^{2+\varepsilon
}dx\right)
\end{equation*}
for $1<q<1+\frac{n+2}{\varepsilon }$.
\end{proposition}

\begin{proof}
From the previous proposition, for $p>\frac{2+n+\varepsilon }{2+n},$ we
proved that $\Psi _{2.\Omega }\in W^{p,k}\left( \Omega ,\parallel x\parallel
^{2+\varepsilon }dx\right) $.

\ 

Therefore, if $\phi $ is a function in $W^{q,k}\left( \Omega ,\parallel
x\parallel ^{2+\varepsilon }dx\right) $ such that $p^{-1}+q^{-1}=1$, we have
the integral estimate :

\ \ 
\begin{equation*}
\dint\limits_{\Omega }\mid \Psi _{2.\Omega }\phi \mid \leq \parallel \Psi
_{2.\Omega }\parallel _{W^{p,k}\left( \Omega ,\parallel x\parallel
^{2+\varepsilon }dx\right) }.\parallel \phi \parallel _{W^{q,k}\left( \Omega
,\parallel x\parallel ^{2+\varepsilon }dx\right) }
\end{equation*}%
where $1<q<1+\frac{n+2}{\varepsilon }$.
\end{proof}

\ \ \ \ 

\begin{proposition}
Over unbounded domain $\Omega \subseteq 
\mathbb{R}
^{n}$, and for $1+\frac{\varepsilon }{n+2}<p<\infty $, with respect to the
weighted measure $d\mu \left( x\right) =\parallel x\parallel ^{2+\varepsilon
}dx$

we have the following norm estimates and norm limit:

\ 
\begin{eqnarray*}
&\parallel &\Psi _{2.\Omega }\parallel _{W^{p,k}\left( \Omega ,\parallel
x\parallel ^{2+\varepsilon }dx\right) }.\parallel \phi \parallel _{W^{\left(
1+\frac{n+2}{\varepsilon }\right) ,k}\left( \Omega ,\parallel x\parallel
^{2+\varepsilon }dx\right) } \\
&\leq &\parallel \Psi _{2.\Omega }\parallel _{W^{\left( 1+\frac{\varepsilon 
}{n+2}\right) ,k}\left( \Omega ,\parallel x\parallel ^{2+\varepsilon
}dx\right) }.\parallel \phi \parallel _{W^{q,k}\left( \Omega ,\parallel
x\parallel ^{2+\varepsilon }dx\right) }
\end{eqnarray*}

and

\begin{eqnarray*}
&&\underset{q\nearrow \left( 1+\frac{n+2}{\varepsilon }\right) }{\lim }%
\left( \parallel \Psi _{2.\Omega }\parallel _{W^{p,k}\left( \Omega
,\parallel x\parallel ^{2+\varepsilon }dx\right) }.\parallel \phi \parallel
_{W^{q,k}\left( \Omega ,\parallel x\parallel ^{2+\varepsilon }dx\right)
}\right) \\
&=&\parallel \Psi _{2.\Omega }\parallel _{W^{p,k}\left( \Omega ,\parallel
x\parallel ^{2+\varepsilon }dx\right) }.\parallel \phi \parallel _{W^{\left(
1+\frac{n+2}{\varepsilon }\right) ,k}\left( \Omega ,\parallel x\parallel
^{2+\varepsilon }dx\right) }
\end{eqnarray*}
\end{proposition}

\begin{proof}
The first part of the proposition follows from the decreasing monotonic
nature of Lebesgue norm with respect to the increase in the index since $%
q\nearrow _{1}^{\left( 1+\frac{n+2}{\varepsilon }\right) }$

and the second follows from the general theory of continuity of Lebesgue
norm.
\end{proof}

\ \ \ \ \ \ 

\begin{corollary}
When $n=2$, we have :

\begin{eqnarray*}
&\parallel &\Psi _{2.\Omega }\parallel _{W^{,pk}\left( \Omega ,\parallel
x\parallel ^{2+\varepsilon }dx\right) }.\parallel \phi \parallel _{W^{\left(
1+\frac{4}{\varepsilon }\right) ,k}\left( \Omega ,\parallel x\parallel
^{2+\varepsilon }dx\right) } \\
&\leq &\parallel \Psi _{2.\Omega }\parallel _{W^{\left( 1+\frac{\varepsilon 
}{4}\right) ,k}\left( \Omega ,\parallel x\parallel ^{2+\varepsilon
}dx\right) }.\parallel \phi \parallel _{W^{q,k}\left( \Omega ,\parallel
x\parallel ^{2+\varepsilon }dx\right) }
\end{eqnarray*}

and 
\begin{eqnarray*}
&&\underset{q\nearrow \left( 1+\frac{4}{\varepsilon }\right) }{\lim }\left(
\parallel \Psi _{2.\Omega }\parallel _{W^{p,k}\left( \Omega ,\parallel
x\parallel ^{2+\varepsilon }dx\right) }.\parallel \phi \parallel
_{W^{q,k}\left( \Omega ,\parallel x\parallel ^{2+\varepsilon }dx\right)
}\right) \\
&=&\parallel \Psi _{2.\Omega }\parallel _{W^{p,k}\left( \Omega ,\parallel
x\parallel ^{2+\varepsilon }dx\right) }.\parallel \phi \parallel _{W^{\left(
1+\frac{4}{\varepsilon }\right) ,k}\left( \Omega ,\parallel x\parallel
^{2+\varepsilon }dx\right) }
\end{eqnarray*}
\end{corollary}

\section{\protect\bigskip\ \ \textbf{Generating kernels: }$\Psi _{l,\Omega
}\left( x\right) $}

\ \ \ In this section, we extrapolate the idea of constructing singular
integral operators as convolutions with fundamental solutions of the Dirac
operator to the once generated by fundamental solutions of higher iterates
of the Dirac operator.

\bigskip\ \ 

Some kernels generate hyper singular integral operators and others form
weaker singular integral operators. It is therefore interesting to look at
differences of these formations from the very constructions of the operators.

\bigskip\ \ 

These functions are constructed by \ recursive ( or iterative ) way from the
fundamental solutions of the Dirac operator and its higher iterates and are
given below:

\begin{equation*}
\Psi _{l,\Omega }\left( x\right) =\left\{ 
\begin{array}{l}
\theta \left( n,l\right) \frac{x}{\omega _{n}\Vert x\Vert ^{n-l+1}},\text{
if }l\text{ is odd} \\ 
\\ 
\frac{\theta \left( n,l\right) }{\omega _{n}\Vert x\Vert ^{n-l+1}},\text{if }%
l\text{ is even \ \ }%
\end{array}%
\right.
\end{equation*}

where $l<n$.

\ 

For a detail study of the consructions of these functions and their
application for constructing complete family of functions and minimal family
of functions, one can see \cite{dr1},\cite{dr2}

\begin{proposition}
For $l<n$ and $\Omega ^{\text{unbdd},\text{ smooth}}\subseteq 
\mathbb{R}
^{n}$, the function $\Psi _{l,\Omega }\in W^{p,k}\left( \Omega _{\varepsilon
},Cl_{n}\right) $ for $\left\{ 
\begin{array}{l}
\frac{n}{n-l}<p<\infty ,\text{ when }l\text{ is odd} \\ 
\\ 
\frac{n}{n+1-l}<p<\infty ,\text{ when }l\text{ is even.}%
\end{array}%
.\right. $

\begin{proof}
For $\Omega $ unbounded and smooth with $\Omega _{\varepsilon }=\Omega
\backslash B\left( x,\varepsilon \right) $ for $\varepsilon >0$\ , using
polar coordinates, the integral 
\begin{equation*}
\dint\limits_{\Omega _{\varepsilon }}\Vert \Psi _{l,\Omega }\left( x\right)
\Vert ^{p}dx
\end{equation*}
is dominated by the integral

\begin{equation*}
C\left( \theta ,n,\omega _{n}\right) \dint\limits_{\varepsilon }^{\infty
}r^{-p\left( n-l\right) +n-1}dr
\end{equation*}%
for $l$ odd with finite integral when the index $p$ satisfies the inequality 
\begin{equation*}
\frac{n}{n-l}<p<\infty
\end{equation*}%
and when $l$ is even, it is dominated by the integral: 
\begin{equation*}
C\left( \theta ,n,\omega _{n}\right) \dint\limits_{\varepsilon }^{\infty
}r^{-p\left( n+1-l\right) +n-1}dr
\end{equation*}%
\ which again is convergent for the indices which satisfy the inequality: 
\begin{equation*}
\frac{n}{n+1-l}<p<\infty
\end{equation*}%
where $C\left( \theta ,n,\omega _{n}\right) $ is some constant that depends
on $n,\theta $ and $\omega _{n}$.
\end{proof}
\end{proposition}

\bigskip\ \ \ \ \ \ 

Thus for $l$ : odd, when we work with this generating kernels, we have the
indices $p$ that depends on $l$ and $n$ and the conjugate index $q$ has the
following limiting values:

as $p\rightarrow \frac{n}{n-l},$ we have : $q\rightarrow \frac{n}{l}$.

\ \ 

Thus, as $\Psi _{l,\Omega }\in W^{p,k}\left( \Omega ,Cl_{n}\right) $\ for $%
\frac{n}{n-l}<p<\infty ,$ the working Sobolev spaces for these kernels are $%
W^{q,k}\left( \Omega ,Cl_{n}\right) $ for $1<q<\frac{n}{l}$ such that $%
p^{-1}+q^{-1}=1$ .

\ 

Therefore for $\phi \in W^{q,k}\left( \Omega ,Cl_{n}\right) $, we have the
convergence of the \textit{sub-singular} or in the literature terminology 
\textit{weakly} singular integral operators :

\begin{equation*}
\dint\limits_{\Omega _{\varepsilon }}\Psi _{l,\Omega }\left( x\right) \phi
\left( x\right) dx
\end{equation*}%
with the usual integral inquality: 
\begin{equation*}
\left( \dint\limits_{\Omega _{\varepsilon }}\Vert \Psi _{l,\Omega }\left(
x\right) \phi \left( x\right) \Vert dx\right) ^{pq}\leq \dint\limits_{\Omega
_{\varepsilon }}\Vert \Psi _{l,\Omega }\left( \ x\right) \Vert
^{p}dx\dint\limits_{\Omega _{\varepsilon }}\Vert \phi \left( \ x\right)
\Vert ^{q}dx
\end{equation*}

\bigskip\ \ 

\ \ For $l$ even, we have the conjugate index $q\rightarrow \frac{n}{l-1}$
as $p\downarrow \frac{n}{n+1-l}$ and since $l<n,$ we have that $\frac{n}{l-1}%
>1$ and therefore, the above inequality holds again.

\ 

Then as convolution, we have :

\begin{proposition}
For $1<q<\frac{n}{l}$,or $1<q<\frac{n}{l-1}$ , the integral operator: 
\begin{equation*}
\dint\limits_{\Omega }\Psi _{l,\Omega }\left( x-y\right) \phi \left(
y\right) dy
\end{equation*}
is a weak-singular integral operator from: 
\begin{equation*}
W^{q,k}\left( \Omega ,Cl_{n}\right) \rightarrow W^{q,k+1}\left( \Omega
,Cl_{n}\right) .
\end{equation*}
\end{proposition}

\begin{proof}
First, as $\Psi _{l,\Omega }\in W^{p,k}\left( \Omega ,Cl_{n}\right) $ for $%
1<p<\infty $, we have that for $\phi \in W^{q,k}\left( \Omega ,Cl_{n}\right) 
$, for $1<q<\frac{n}{l}\left( \text{ for }l\text{ odd}\right) $ or for $1<q<%
\frac{n}{l-1}\left( \text{for }l\text{ even}\right) $ with $p^{-1}+q^{-1}=1$

such that the integral 
\begin{equation*}
\dint\limits_{\Omega _{\varepsilon }}\Psi _{l,\Omega }\left( x-y\right) \phi
\left( y\right) dy
\end{equation*}
is convergent but singular with out the puncture .

\ \ 

The convolution is the usual Teodorescu transform which has the mapping
property :

\ \ 
\begin{equation*}
\Psi _{l,\Omega }\ast \phi :W^{q,k}\left( \Omega ,Cl_{n}\right) \rightarrow
W^{q,k+1}\left( \Omega ,Cl_{n}\right) .
\end{equation*}
\end{proof}

\ \ 

\begin{proposition}
In the usual $3-D$ Euclidean space, if $l=3$, then we can not work on the
usual generalized Hilbert space $W^{2,k}\left( \Omega ,Cl_{n}\right) $.
\end{proposition}

\begin{proof}
For such a setting, we have that $3<p<\infty $ and therefore the working
function spaces will have conjugate Sobolev indices with range $1<q<\frac{3}{%
2}$, in which the index $2$ is not included.

\ \ \ 

Therefore the Sobolev space of index $2$ which is the generalized Hilbert
space $W^{2,k}\left( \Omega ,Cl_{n}\right) $ is no more a viable space.
\end{proof}

\ \

\end{document}